\documentclass[12pt]{article}
\usepackage{amssymb}
\def\be{\begin{equation}}
\def\ee{\end{equation}}

\def\R{{\sf I\kern-.2em R}}
\def\N{{\sf I\kern-.2em N}}
\def\C{\kern.1em{\raise.47ex\hbox{$\scriptscriptstyle
$}}\kern-.40em{\sf C}}
\def\Z{{\sf Z\kern-.32em Z}}

\def\be{\begin{equation}}
\def\ee{\end{equation}}

\newtheorem{theorem}{\noindent Theorem}
\newtheorem{lemma}{\noindent Lemma}
\newtheorem{definition}{\noindent Definition}
\newtheorem{corollary}{\noindent Corollary}

\newtheorem{statement}{\noindent Statement}
\newtheorem{remark}{\noindent Remark}

\begin{document}

\begin{titlepage}

\begin{center}
{\LARGE \bf RANDOM AND UNIVERSAL METRIC SPACES}
\bigskip
\bigskip

{\large A.~M.~VERSHIK $^{1}$}
\bigskip

{\sl
$^{1}$ Steklov Institute of Mathematics at St.~Petersburg,\\
Fontanka 27,
119011 St.~Petersburg, Russia.\\

Partially supported by RFBR, grant 02-01-00093, and by CDRF, grant RMI-2244).}
\end{center}
\bigskip
\bigskip
\begin{center} ABSTRACT
\end{center}

   We introduce a model of the set of all Polish (=separable complete metric)
   spaces: the cone $\cal R$ of distance matrices, and consider  geometric
   and probabilistic problems connected with this object. The notion of the
   universal distance matrix is defined and we proved that the set of
   such matrices is everywhere dense $G_{\delta}$ set in weak topology
   in the cone $\cal R$.
   Universality of distance matrix is the necessary and sufficient condition
   on the distance matrix of the countable everywhere dense set of so called
   universal Urysohn space which he had defined in 1924 in his last paper.
   This means that Urysohn space is generic in the set of all Polish spaces.
   Then we consider metric spaces with measures (metric triples)
   and define a complete invariant: its - matrix distribution.
   We give an intrinsic characterization of the set of matrix distributions,
   and using the ergodic theorem, give a new proof of Gromov's ``reconstruction
   theorem'. A natural construction of a wide class of measures on the
   cone $\cal R$ is given and for these we show that {\it with probability
   one a random Polish space is again the Urysohn space}. There is a
   close connection between these questions, metric classification of
   measurable functions of several arguments,
   and classification of the actions of the infinite symmetric group
   (\cite{V1,V2}).

\begin{center} CONTENT
\end{center}

1.The cone of the distance matrices as the set of all Polish spaces.

2.Geometry and topology of the cone of distance matrices.

3.Universal matrices and Urysohn space.

4.Matrix distribution as complete invariant of the metric triples
  and its characterization.

5.General classification of the measures of the cone
  of the distance matrices, examples.

Bibliography.

\end{titlepage}
\newpage

\section{Introduction: The cone $\cal R$ of distance matrices
as a set of all the Polish spaces}

Consider the set of all infinite real matrices
 $${\cal R} =\{\{r_{i,j}\}_{i,j=1}^{\infty}:  r_{i,i}=0,\,r_{i,j}\geq 0,\,
 r_{i,j}=r_{j,i},\,  r_{i,k}+r_{k,j}\geq r_{i,j},{\rm \,\,for\,\,} i,j,k=1,2, \dots\}$$

We will call the elements of $\cal R$  {\it distance matrices}.
Each such matrix defines a semi-metric on the set of natural numbers
$\bf N$. We allow zeros away from the principal diagonal, so
in general  $\rho$ is  only a semi-metric.
If matrix has no zeros away from the principal diagonal we will call
it {\it a proper distance matrix}.

The set of all distance matrices is a weakly closed convex cone in the
real linear space $Mat_{\bf N}({\bf R})={\bf R^N}^2$  endowed with the
ordinary weak topology.
We always consider the cone ${\cal R}$ with this topology and will call it the
cone of distance matrices. The subset of proper distance matrices
is everywhere dense open subcone in the cone ${\cal R}$.

If the distance matrix $r$ is proper then the completion of the
metric space $({\bf N},r)$ is a complete separable metric space
(=Polish space) ($X_r,\rho_r)$ with a distinguished everywhere dense
countable subset
$\{x_i\}$ which is the image of the natural numbers in the completion.
 A general distance matrix
(with possible zeros away from the diagonal)
defines on the set of natural numbers structure of
{\it semi-metric space}. By the completion of $({\bf N},r)$ in this case
we mean the completion of the corresponding quotient metric
space of the classes
of points with zero distances. For example the zero matrix
is a distance matrix
on the natural numbers  with zero distances between each two numbers
and its ``completion'' is the singleton metric space. Thus finite metric
spaces also could be considered in this setting.

Suppose now that we have some  Polish  space $(X, \rho)$,
equipped with the orderes everywhere dense countable set
$\{x_i\}_{i=1}^\infty$.
Defining the matrix $r=\{r_{i,j}\} \in \cal R$ by
$r_{i,j}=\rho(x_j,x_j), i,j =1\dots$  we obtain a proper distance matrix.
which we interpret as a {\it metric on the set of natural numbers}.
Clearly this distance matrix analogously to structural constant
in the algebraic situation, contains complete information about
the original
space $(X,\rho)$ because  $(X,\rho)$ is the canonical
completion of the set of naturals with this metric.
Any invariant property of the metric space
(topological and homological etc.) could be expressed in terms of the
distance matrix for any dense countable subset of that space.
We will study the theory of Polish spaces from this point of view and
consider the cone of distance matrices ${\cal R}$ as {\it the universe of
all separable complete metric spaces} with a fixed dense countable subset
and study the properties of the metric spaces as well as properties
of whole set of its using thuis cone.
We can view $\cal R$ as a ``fibering'', whose base is
the collection of all individual Polish spaces, and the fiber over a given space is
the set of all countable ordered dense subsets in this metric
space. Because of the universality of the Urysohn space $\cal U$ (see below)
the set of all closed subsets of $\cal U$  could be considered as
a base of this bundle.
Thus the space $\cal R$  plays the role of a "tautology fibration" over
the space of classes of isometric Polish spaces, analogous to common
topological constructions of a tautology fiber bundles.

The question arises: what kind of distance matrix is ``generic''
in the sense of the topology of $\cal R$. One of the main results
(section 3) is the Theorem 1, which is a generalization of Urysohn's
results and which asserts that Urysohn space is generic
(=dense $G_{\delta}$ set in $\cal R$)
in this sense. The main tool is the notion of {\it universal distance
matrices}, an example of such matrix was used in indirect way
by P.S.Urysohn in his pioneer paper \cite{U} for the proof of
exitence of the universal metric space. An explicit formulation
of the notion of universality of the distance matrix is given
in Statement 1 (section 3.1).
We give a new version of his main results and a new proof in the section 3.
Related consideration of the Urysohn space can be found also in the papers
\cite {Ka,Bo,Us,Pe}. I want to point out that the fact that during
almost 70 years Urysohn's paper \cite{U} with this result
was out of attention of the mathematicians is astonishimg;
I do not know any text-book or monographs on general topology
in which Urysohn universal metric space was mentioned!

Introduce a partition $\xi$ of the $\cal R$ into the equivalence
classes of distance matrices with isometric completions.
The quotient space over the partition $\xi$ (or space of the classes
of equivalence)
is the set of the {\it isometry classes} of the Polish spaces.
As was conjectured in  \cite{V1} and proved in the paper \cite{CK01}
the quotient by this equivalence relation is not "smooth", in
the sense that it has no good Borel or topological structure and thus
the problem of the classification of the Polish spaces up to
isometry is ``wild''.  At the same time the restricted
problem for the case of compact Polish spaces is smooth (see \cite{G})
and the space of all isometry classes of compact metric spaces
has a natural topology. Surprisingly, if we consider the problem
of classification of the Polish spaces with measure (metric triples)
up-to isometry which preserves the measures, this classification
is ``smooth'', and we will consider in detales a complete invariant
(``matrix distribution'') of metric triples (section 4).
One direction -the completeness of this invariant - was proved
in the book by M.Gromov \cite{G}; we will give another proof
of his reconstruction theorem based on ergodic methods and
a new description of the invariant. Then we prove a theorem about
the precise description of the matrix distribution of the metric
triples (section 4) as a measure on $\cal R$.
 The section 2 is devoted to the elementary geometry of the
space $\cal R$ which we use throughout  all the paper, and
especially in the section 5 in which we consider the various types
of measures on the cone $\cal R$, and methods of the constructing
of them. The measure on the above cone is nothing more than
a random metric on the set of naturals numbers.
Thus we can construct a ``random'' metric space as the result
of completion of the random metric on the natural numbers.
In this way we prove that loosely speaking, a Polish space
randomly constructed, by a very natural procedure gives us
with probability 1 again universal Urysohn space.
We can say, that {\it the random space is universal space}
(see \cite{V5}).

One of the previous simple analogy of such theorems is the theorem due to
P.Erd{\"o}s and A.R\'{e}nyi about random graphs (see \cite{ER,Ca}).
The results of the paper about the  genericity of the Urysohn space
(Theorem 1) and probabilistic typicalness of its (Theorem 7) show
that these two properties  coincide in the category of the Polish spaces
as well as in more simple case of the graphs. Perhaps,
this coincidence has more a general and deep feature and
takes place in the other categories.
As a similar facts recall universality of Poulsen simplex (\cite{LO})
and of the Guraij's Banach space (\cite{Lu}) ( Y.Beniamini's remark),
exsitence of the group which is universal in the class of finite groups
homogeneous Hall's group.

Many questions about Urysohn space remain open, it is not clear if it
is contractable or not there are no good realization of it; one of the
main question is to construct a natural probability measures in the
space. The group of the isometries of Urysohn space is also
very intriguing object (see \cite{Us,Pe}). We will discuss these
questions elsewhere.                	

\newpage
\section{Geometry and topology of the cone ${\cal R}$}
\subsection{Convex structure}

Analogously to $\cal R$ let us define the finite dimensional
cone ${\cal R}_n$ of distance matrices
of order $n$. Cone ${\cal R}_n$ is a polyhedral cone
inside the positive orthant in $ Mat_n({\bf R})\equiv {\bf R}^{n^2}$.
Denote by ${\bf M}^s_n({\bf R})\equiv {\bf M}^s_n$ of the space
{\it symmetric matrices with zeros on the  principal diagonal}.  We have
${\cal R}_n \subset {\bf M}^s_n$ and the latter space is evidently
the linear hull of the cone: span$({\cal R}_n)={\bf M}^s_n$,
because the interior of ${\cal R}_n$
is not empty. It is clear that $span({\cal R})\subset {\bf M}^s_{\bf N}$,
where  ${\bf M}^s_{\bf N}$ is the space of all real infinite
symmetric matrices with zero principal diagonal, the geometry
of the cone $\cal R$ is very complicated.

Each  matrix $r \in {\cal R}_n$ defines a (semi)metric
space $X_r$ on  the  of $n$-point set.

Define the projection $$p_{m,n}:{\bf M}^s_m \longrightarrow {\bf M}^s_n, m >n $$
which associates with the matrix $r$ of order $m$ its NW-corner of order $n$.
The  cones ${\cal R}_n$ are consistent with the projections i.e.
$p_{n.m}: p_{m,n}({\cal R}_m)={\cal R}_n$.
The projections $p_{n,m}$  extend naturally to the space of infinite
symmetric matrices with zero diagonal - $p_n: M^s_{\bf N}  \longrightarrow M^s_n({\bf R})$,
and $p_n$ also preserve the cones: $p_n({\cal R})={\cal R}_n$.
It is clear that ${\cal R}$ is the inverse limit as topological space (in weak topology)
of  the system $({\cal R}_n,\{p_n\})$.

An important but evident property of the cone ${\cal R}_n$
is its invariance under the action of the symmetric group  $S_n$
simultaneously permuting the rows and columns of the
matrices.

Let us consider the geometrical structure of ${\cal R}_n$ and ${\cal R}$.

For the first two dimensions we have ${\cal R}_1=\{0\}),
{\cal R}_2={\bf R}$.
It is interesting to describe the extremal rays (in the sense
of convex geometry) of the convex polyhedral cone ${\cal R}_n, n=3, \dots, \infty$.  This is a well-known problem - see \cite{DL, Av} and
the list of literature there.
Each extremal ray in ${\cal R}_n, n \leq 4$ is of the type
$\{\lambda \cdot l: \lambda \geq 0\}$,
where $l$ is a symmetric $0-1$-
distance matrix which corresponds to the semi-metric space
whose quotient metric space has just two points.
For $n \geq 5$ there are extremal rays of other type.
The complete description of the set extremal rays is rather a difficult
and very interesting combinatorial problem. The most important question for us
concerns to the asymptotic properties of cone ${\cal R}_n$ and especially the
description of the set of extremal rays of the infinite dimensional
cone ${\cal R}$. It happens that this set is a  dense $G_{\delta}$ in
${\cal R}$ and some of the so called universal distance matrices (see par 3.)
are extremal. This is in consistent with the estimation in \cite{Av}
of the number of extremal rays of ${\cal R}_n$ which grows very rapidly.
The algebro-geometric structure and stratification
of the cones ${\cal R}_n$  as semi-algebraic sets.
are also very intriguing.
In order to clarify topological and convex structure
of the cones ${\cal R}_n$
we will use an {\it inductive description} of these cones and
will study it in the next subsection.
	
\subsection{Admissible vectors and structure of the $\cal R$}

Suppose $r=\{r_{i,j}\}_1^n$ is a distance matrix of order $n$
($r \in {\cal R}_n$),
choose a vector $a \equiv \{a_i\}_{i=1}^n \in {\bf R}^n$
such that if we attaching to the matrix $r$ with vector $a$ as the last
column and the last row then the new matrix of order $n+1$ still belongs
to  ${\cal R}_{n+1}$. We will call such vector {\it an admissible vector}
for fixed distance matrix $r$ and denote the set of of all admissible
vectors for $r$ as $A(r)$. For given $a \in A(r)$ denote as $(r^a)$,
distance matrix of order $n+1$ obtained from  matrix $r$
adding vector $a \in A(r)$ as the last row and column.
It is clear that $p_n(r^a)=r$. The matrix $r^a$ has the form

 $$r^a = \left( \begin{array}{ccccc}
 0    &r_{1,2}& \ldots&  r_{1,n}&a_1\\
 r_{1,2} &0& \ldots &r_{2,n} &a_2 \\
 \vdots &  \vdots&
\ddots & \vdots &\vdots\\
r_{1,n}&r_{2,n}& \ldots &0  &a_n\\
  a_1&  a_2  & \ldots & a_n& 0
\end{array}\right)$$

The (semi)metric space $X_{r^a}$ corresponding to matrix $r^a$
is an extension
of $X_r$: we add one new point $x_{n+1}$ and $a_i, i=1 \dots n$
is the distance
between $x_{n+1}$ and $x_i$.
The admissibility  of $a$ is equivalent to the following set
of inequalities:
the vector $a=\{a_i\}_{i=1}^n$ must satisfy to the series of
triangle inequalities
for all $i,j,k=1,2 \dots n$; (matrix $\{r_{i,j}\}_{i,j=1}^n$ is fixed):
\begin{equation}
\label{trivial} |a_i-a_j|\leq r_{i,j}\leq a_i+a_j
\end{equation}

So for given distance matrix $r$ of order $n$ the set of
admissible vectors
is $A(r)=\{\{a_i\}_{i=1}^n:|a_i-a_j|\leq r_{i,j}\leq a_i+a_j, i,j=1 \dots n\}$.
It makes sense to mention that we can view on the vector
$a=\{a_i\}$ as a {\it Lipshitz function} $f(.)$
on the space $X_r=\{1,2 \dots n\}$ with $r$ as a metric:
$f_a(i)=a_i$  with Lipshitz constant equal $1$. This point of view
helps to consider a general procedure of extension of metric space.

Geometrically the set $A(r)$ can be identified with
the intersection of cone ${\cal R}_{n+1}$ and the affine subspace
which consists of matrices of order $n+1$ with
given matrix $r$ as the NW-corner of order $n$.
It is clear from the linearity of inequalities that
the set $A(r)$ is an
unbounded closed convex polyhedron in ${\bf R}^n$.
If $r_{i,j} \equiv 0, i,j=1 \dots n\geq 1$,
then $A(r)$ is diagonal: $A(0)=\Delta_n \equiv \{(\lambda,
\dots \lambda): \lambda \geq 0\}\subset {\bf R}^n_+$.
Let us describe the structure of $A(r)$ more carefully.

\begin{lemma}
{For each proper distance matrix $r$ of order $n$
the set of admissible vectors $A(r)$
is a closed convex polyhedron in the orthant ${\bf R}_+^n$,
namely it is a Minkowski sum:
$$A(r)= M_r +\Delta_n,$$
where  $\Delta_n $ is the half-line of constant vectors in the
space ${\bf R}_+^n$,
 and $M_r$ is a compact convex polytope of dimension
$n$. This polytope $M_r$ is the convex hull of extremal points
of the polyhedron $A(r)$:   $M_r= conv(ext A(r))$.}
\end{lemma}

\begin{proof}
{The set $A(r) \subset {\bf R}^n$ is the intersection of finitely many
closed
half-spaces, and evidently it does not contain straight lines, so, by a
general theorem of convex geometry $A(r)$ is a sum of the convex
closed polytope
and some cone (which does not contain straight lines) with the vertex
at origin.
This convex polytope is the convex hull of the extremal points
of convex set $A(r)$.
But this cone is half-line of the constant (diagonal) nonnegative
vectors in ${\bf R}^n$ because if it contains
half-line which is different from the diagonal then
the triangle inequality is violated. The dimension of $A(r)$ equal to $n$
for proper distance matrix; in general it
depends on the matrix $r$ and could be less than $n$ for some matrix $r$;
while the dimension of $M_r$ is equal to $\dim A(r)$ or to $\dim A(r)-1$.
The assertion about topological structure
of $A(r)$ follows from what was claimed above.}
\end{proof}

The next lemma asserts that this correspondence
$r \rightarrow A(r)$ is covariant under the action of symmetric
group in ${\bf R}^n$.  The proof is evident.

\begin{lemma}
{ For any $r\in {\cal R}_n$ we have coincideness of the sets:
$A(grg^{-1})=g(A(r))$, where $g \in S_n$ is element of symmetric
group $S_n$ which acts in a natural way on the space of matrices
$M_{\bf N}(\bf R)$ as well as on the space of the convex subsets
of the vector space ${\bf R}^n$.}
\end{lemma}

The convex structure of polytopes $M_r, A(r)$ is very interesting
and seems to have not been studied before. For dimensions higher
than 3 {\it combinatorial type}
of the polytope $M_r$ hardly depends on $r$ but for
dimension three the combinatorial type of polytopes $M_r$,
and consequently the
combinatorial structure of the sets $A(r)$ is the same for
all proper distances matrices $r$.

{\bf Example}
For $n=3$ the description of the set $A(r)$ and of its extremal
points is the following.
Let $r$ be a matrix
 $$r = \left( \begin{array}{ccc}
  0 & r_{1,2} & r_{1,3}\\
  r_{1,2}& 0 & r_{2,3} \\
  r_{1,3} &r_{2,3}&0 \\
\end{array}\right)$$
Denote $r_{1,2}=\alpha,  r_{1,3}= \beta,  r_{2,3}=\gamma$, then
 $r^a$:

 $$r^a = \left( \begin{array}{cccc}
 0    &\alpha  &  \beta & a_1\\
 \alpha &0& \gamma & a_2 \\
 \beta & \gamma &0  & a_3\\
  a_1&  a_2  &  a_3 &  0
\end{array}\right)$$

Denote $\delta=\frac{1}{2}(\alpha+\beta+\gamma)$
There are seven extremal points $a=(a_1,a_2,a_3)$ of $A(r)$ :
the first one is a vertex which is the closest to the origin:
$(\delta-\gamma, \delta-\beta, \delta-\alpha)$,
then another three non degenerated extremal points:
$(\delta,\delta-\alpha, \delta-\gamma),
(\delta-\beta, \delta, \delta-\alpha),(\delta-\gamma, \delta-\beta, \delta)$,
and  three degenerated extremal points
$ (0, \alpha, \beta), (\alpha, 0, \gamma), (\beta, \gamma, 0)$.

If $\alpha=\beta=\gamma=1$
then those seven points are as follows

$(1/2,1/2,1/2), (3/2,1/2,1/2),(1/2,3/2,1/2),(1/2,1/2,3/2),
(0,1,1),(1,0,1),(1,1,0).$

Remark that all non-degenerated extremal points in the example
defines the finite metric spaces which can not be isometrically
embedded to Euclidean space.

\subsection{Projections and isomorphisms}

Let $r$ be a distance matrix of order $N$ and $p_n(r)$ its NW-corner of order $n<N$.
Then we can define a projection
$\chi^r_n$ of $A(r)$ to $A(p_n(r))$:
$\chi^r_n:(b_1, \dots b_n, b_{n+1},\dots b_N) \mapsto (b_1 \dots b_n)$.
(We omit the index $N$ in the notation  $\chi^r_n$).
The next simple lemma plays a very important role for our construction.
\begin{lemma}
{Let  $r\in {\cal R}_n$ be a distance matrix of order $n$. For any two vectors
$a=(a_1, \dots a_n) \in A(r)$ and  $b=(b_1, \dots b_n) \in A(r)$ there exists
a real nonnegative number $h \in {\bf R}$ such that vector
${\bar b}=(b_1,\dots b_n, h) \in A(r^a)$
(and also ${\bar a}=(a_1, \dots a_n, h) \in A(r^b)$).}
\end{lemma}

\begin{corollary}
{For each $r \in {\cal R}_n$ and $a \in A(r)$ the map
$\chi^r_{n+1,n}: (b_1, \dots b_n, b_{n+1}) \mapsto (b_1 \dots b_n)$ of  $A(r^a) \to A(r)$
is the epimorphism of $A(r^a)$ onto $A(r)$ (by definition $p_{n+1,n}(r^a)=r$)}
\end{corollary}

\begin{proof}
{The assertion of this lemma as we will see, follows from  simple geometrical observation:
suppose we have two finite metric spaces
$X=\{x_1, \dots x_{n-1}, x_n\}$ with metric $\rho_1$
and $Y=\{y_1, \dots y_{n-1},y_n \}$ with metric $\rho_2$. Suppose the subspaces
of the first $n-1$ points $\{x_1, \dots x_{n-1}\}$ and $\{y_1, \dots y_{n-1}\}$
are isometric, i.e. $\rho_1(x_i,x_j)=\rho_2(y_i,y_j), i,j=1, \dots n-1$.
Then there exists a third space $Z=\{z_1,\dots z_{n-1},z_n,z_{n+1}\}$ with metric
$\rho$ and two isometries  $I_1,I_2$ of both spaces $X$ and $Y$ to the space
$Z$ so that $I_1(x_i)=z_i,I_2(y_i)=z_i, i=1, \dots n-1, I_1(x_n)=z_n,
I_2(y_n)=z_{n+1}$.
In order to prove existence of $Z$ we need to show that it is possible to
define only nonnegative number $h$ which will be
the distance $\rho(z_n,z_{n+1})=h$ between $z_n$ and $z_{n+1}$ (images of $x_n$ and $y_n$
in $Z$ correspondingly)
such that all triangle inequalities
are valid in the space $Z$.
The existence of $h$ follows from the inequalities:
$$\rho_1(x_i,x_n)-\rho_2(y_i,y_n) \leq \rho_1(x_i,x_j)+\rho_1(x_j,x_n)-\rho_2(y_i,y_n)=$$
$$=\rho_1(x_j,x_n)+\rho_2(y_i,y_j)-\rho_2(y_i,y_n)\leq \rho_1(x_j,x_n)+
\rho_2(y_j,y_n)$$
for all $i,j=1, \dots n-1$.

 Consequently $$\max_i |\rho_1(x_i,x_n)-\rho_2(y_i,y_n)| \equiv M
\leq  m \equiv  \min_j (\rho_1(x_j,x_n)+\rho_2(y_j,y_n)).$$
Thus, a number $h$ could be chosen as an arbitrary number from
the nonempty closed interval $[M,m]$ and we set
$ \rho(z_n,z_{n+1}) \equiv h$;
it follows from the definitions that all triangle inequalities
are satisfied.
Now suppose we have a distance matrix $r$ of order $n-1$ and
an admissible vector $a \in A(r)$,
so we have a metric space $\{x_1, \dots x_{n-1}, x_n\}$
(the first $n-1$ points correspond
to the matrix $r$, and whole space - to the extended matrix $r^a$.
Now suppose we choose
another admissible vector $b \in A(r)$, and giving distance
matrix $r^b$ of the space
$\{y_1, \dots y_{n-1},y_n\}$, where the subset of first $n-1$
points is isometric
to the space $\{x_1, \dots x_{n-1}\}$. As we proved we can define
space $Z$ whose
distance metric $\bar r$ of order $n+1$ gives the required property.}
\end{proof}

Now we can formulate a general assertion about the projections $\chi^r$.
\begin{lemma}
{For arbitrary natural numbers $N$ and $n$, $N>n$,
and any $r \in {\cal R}_N$
the map $\chi^r_n$ is epimorphism of $A(r)$ onto $A(p_n(r))$.
In other words for each $a=(a_1, \dots a_n) \in A(p_n(r))$
there exists a vector $(b_{n+1},\dots b_N)$ such that
 $b=(a_1,\dots a_n,  b_{n+1}, \dots b_N) \in A(r)$.}
\end{lemma}
\begin{proof}
{The above proof shows how to define the first number $b_{n+1}$.
But the projection $\chi^r_n$ seen as a map from  $A(r),
r\in {\cal R}_N$ to $A(p_n(r))$
is the product of projections $\chi^r_n\cdots \chi^r_{N-1}$.
Because each factor is epimorphism the product is epimorphism also.}
\end{proof}

It is convenient for our goals to represent the infinite distance matrix
$r \equiv \{r_{i,j}\} \in {\cal R}$
as a sequence of admissible vectors of increasing lengths:
\begin{equation}
r(1)=\{r_{1,2}\},
r(2)=\{r_{1,3},r_{2,3}\}, \dots r(k)=\{r_{1,k+1},r_{2,k+1}, \dots r_{k,k+1}\}\dots,
k=1,2 \dots,
\end{equation}
satisfying conditions
$r(k) \in A(p_k(r))$, (recall that $p_k(r)$ is the NW-projection
of matrix $r$ on the space
${\bf M}^s_k$ defined above), so each vector $r(k)$ is admissible
for the {\it previous distance matrix $p_k(r)$}.
We can consider the following sequence of the cones and maps:

\begin{equation}
\label{trivial}0={\cal R}_1\stackrel{p_2}{\longleftarrow}
{\cal R}_2={\bf R}_+\stackrel{p_3}{\longleftarrow}{\cal R}_3
{\longleftarrow} \dots
{\longleftarrow}{\cal R}_{n-1}\stackrel{p_n}{\longleftarrow}{\cal R}_n{\longleftarrow} \dots
\end{equation}
the projection $p_n$ here is the restriction of the projection defined
above onto the cone ${\cal R}_n$.
The  preimage of the point $r \in {\cal R}_{n-1}$ (fiber over $r$)
is the set $A(r)$ which described in Lemma 1.
Note  that this is not fibration in the usual sense:
the preimages of the points could
even not be homeomorphic to each other for  various $r$
(even dimensions could be different). But that sequence defines
allows to define cone $\cal R$ as an inverse limit of
the cones ${\cal R}_n$.
We will use the sequence (3) in order to define the measures
on the cone ${\cal R}$ in the spirit of the theory of
the Markov processes.

\newpage
\section{Universality and Urysohn space}
\subsection{Universal distance matrices}

The following definition plays a crucial role.
\begin{definition}
{1.An infinite proper distance matrix
$r=\{r_{i,j}\}_{i,j=1}^{\infty} \in {\cal R}$
is called a \uchyph=0 UNIVERSAL distance matrix if the following condition is satisfied:

for any $\epsilon >0$, $n \in \bf N$ and for any
vector $a=\{a_i\}_{i=1}^n \in A(p_n(r))$
there exists $m \in \bf N$  such that
$\max_{i=1 \dots n}|r_{i,m}-a_i| < \epsilon.$

In another words: for each $n \in \bf N$ the set of vectors
$\{\{r_{i,j}\}_{i=1}^n\}_{j=n+1}^{\infty}$ is everywhere dense in the set
of admissible vectors  $A(p_n(r))$.

2.An infinite proper distance matrix
$r=\{r_{i,j}\}_{i,j=1}^{\infty} \in {\cal R}$
is called a weakly universal distance matrix if
for any $n \in \bf N$ the set of all submatrices
$\{r_{i_k,i_s}\}_{k,s=1}^n$ of the matrix  $r$ of order $n$
over all $n$-tuples $\{i_k\}_{k=1}^n \subset \bf N$ is dense
in the cone ${\cal R}_n$.}
\end{definition}

Let us denote the set of universal distance  matrices by $\cal M$.
We will prove that $\cal M$ is not empty but before we
formulate some properties of universal matrices.

\begin{lemma}
{Each universal distance matrix is weakly universal.
There exist nonuniversal but weakly universal distance matrices.
}
\end{lemma}

\begin{proof}
{Choose any distance matrix $q \in {\cal R}_n$; we will prove that for given positive
$\epsilon$ it is possible to find a set $\{i_k\}_{k=1}^n \subset \bf N$
such that $\max_{k,s=1,\dots n}|r_{i_k,i_s} - q_{k,s}|<\epsilon$.
Because $r^1=\{r_{1,1}=0\}$ then $A(r^1)={\bf R}_+$
(see 2.2), and by universality of $r$ the sequence
$\{r_{1,n}\}_{n=2}^\infty$ must be dense in ${\bf R}_+ $, so we can
choose some
$i_1$ such that $|r_{1,i_1}-q_{1,2}|<\epsilon$, then using density
of the columns of length 2 which follows from the universality condition
we can choose a natural number $i_2$ such that
$|r_{1,i_2}-q_{1,3}|< \epsilon, |r_{2,i_2}-q_{2,3}|< \epsilon$, etc.

There are many examples of weakly universal but nonuniversal
distance matrices. The distance matrix of the arbitrary countable everywhere
dense set of any universal but not homogeneously universal (see below)
Polish spaces (like $C([0,1])$) gives such a counterexample,
but the simplest one is the distance matrix of the disjoint union of all
finite metric spaces with the rational distance matrices (B.Weiss's example).}
\end{proof}

The following corollary of universality gives useful tool for
tre studying of the Urysohn's space:

\begin{corollary}($\epsilon$-extension of the isometry)
{Suppose $r$ is an infinite universal distance matrix and $q$ is a finite
distance matrix of order $N$ such that for some $n<N,
r_{i,j}=q_{i,j}, i,j=1 \dots n$. Denote $i_k=k, k= 1 \dots n$.
 Then for any positive $\epsilon$ there exist the natural numbers
$i_{n+1} \dots i_N$ such that
$\max_{k,s=1 \dots N}|r_{i_k,i_s}-q_{k,s}|<\epsilon$.
In another words, we can enlarge the set of the first $n$
natural numbers with
some set of $N-n$ numbers: $i_{n+1}, \dots, i_N$ in such a way
that the distance matrix of whole set $i_1=1, i_2=2 \dots i_n=n,
i_{n+1}, \dots i_N$ differs with the distance matrix $q$ (in norm)
less than $\epsilon$.

Conversely, if the infinite distance matrix $r$ has property above
for any finite distance
submatrix $q$ and for any positive $\epsilon$ then matrix $r$
is a universal matrix.}
\end{corollary}

\begin{proof}
{For $N=n+1$ the claim follows directly from the definition of universality of $r$,
then we can use induction on $N$. The second claim follows from the definition.}
\end{proof}

Once more reformulation of the notion of universality
uses the term of group action. Suppose $q \in {\cal R}_n$; denote
${\cal R}^n(q)$ the set of all $r \in \cal R$, with NW-corner
equal to the matrix $q$. Consider the group $S^n_{\infty}$ of
permutations, which preserve as fixed the first $n$ rows and columns
of the matrices from $\cal R$ and consequently map the
set ${\cal R}^n(q)$ to itself.
The following criteria of the universality is a direct corollary of
the definition :
\begin{statement}
{A matrix $r \in \cal R$ is universal iff for each $n$
the orbit of $r$ under the action of the group $S^n_{\infty}$
is everywhere dense in ${\cal R}^n(r^n)$ in weak topology,
here $r^n$ is the NW-corner of matrix $r$ of the order $n$.
The matrix $r$ is weakly universal iff its  $S_{\infty}$-orbit
is everywhere dense in $\cal R$.}
\end{statement}

From other side the existence of universal distance matrix as well as existence
of Urysohn space is not evident. We simplify and a little bit
strengthen Urysohn's existence theorem and prove the following

\begin{theorem}
{The set $\cal M$ of the universal matrices is nonempty.
Moreover, this set is  everywhere dense $G_{\delta}$-set in the cone $\cal R$
in the weak topology.}
\end{theorem}

\begin{proof}
{We will use the representation described in the lemmas in previous
section for construction
of at least one universal proper distance matrix in the cone ${\cal R}$.

Let us fix sequence  $\{m_n\}_{n=1}^\infty$ of natural numbers in which each
natural number occurs infinitely many times and for each $n,  m_n \leq n; m_1=1$.
For each proper finite distance matrix $r \in {\cal R}_n$ let us choose an ordered countable dense
subset $\Gamma_r \subset A(r)$ of the vectors with positive coordinates:
$\Gamma_r=\{\gamma^r_k\}_{k=1}^\infty \subset A(r){\subset \bf R}^n$ and
choose any metric in  $ A(r)$, say the Euclidean metric.

The first step consists of the choice of positive real number
$\gamma_1^1 \in \Gamma_1 \subset A(0)={\bf R}_+^1$
so that we define a distance matrix $r$ of order $2$ with element $r_{1,2}=\gamma_1^1$.

Our construction of the universal matrix $r$ is inductive one,
it used the representation of matrix as a sequence of admissible
vectors $\{r(1),r(2),\dots\}$ with increasing lengths (formula (2)).
The conditions on the vectors are as follows $r(k) \in A(p_k (r_{k+1}))$.
The sequence of the corresponding matrices $r_n, n=1 \dots $
stabilizes to the infinite matrix $r$.
Suppose  after  $(n-1)$-th step we obtain a finite matrix $r_{n-1}$;
then we choose a new admissible vector $r(n) \in  A(r_{n-1})$. The choice of
this vector (denote it $a$) is defined by the condition that
the distance (in norm) between the projection $\chi^r_{m_n}(a)$
of the vector $a$ onto the subspace of admissible vectors $A(r_{m_n})$
and the point $\gamma^{m_n}_s \in \Gamma_{r_{m_n}} \subset A(r_{m_n}) $
must be less than $2^{-n}$, where $s =|i: m_i=m_n, 1\leq i \leq n|+1$:
$$\|\chi^r_{m_n}(a) - \gamma^{r_{m_n}}_s\|<2^{-n}.$$
 Recall now  that the projection $\chi^r_{m_n}$ is an epimorphism
from $A(r)$ to $A(p_{m_n}(r))$,
(Lemma 4), hence a vector $a \in A(r_n)$ with these properties does exist.
The number $s$ is just the number of the points of $\Gamma_{r_{m_n}}$
 which occur on the previous steps of the construction.
After infinitely many steps we obtain the infinite distance
matrix $r$.

 Universality of $r$ is evident, because for each $n$ projection $\chi^r_n$
 of vectors $r(k), k=n+1 \dots$  is a dense set in $A(r_n)$ by construction.
 This proves the existence of the universal matrix.

Now notice that the property of universality of the matrix are preserved
under the action of finite simultaneous permutations of the rows
and columns, and also under the NW-shift which cancels the first
row and first column of the matrices.
Also the set of universal matrices $\cal M$ is invariant under the changing
of the finite part of the matrix. Consequently, $\cal M$ contains together with
the given matrix also its permutations and shifts. But because of the weak
universality of any universal distance matrices $r$,
even the orbit of matrix $r$ under the action of the group of permutations
$S_{\bf N}$ is everywhere dense in $\cal R$ in weak topology. A fortiori
$\cal M$ is everywhere dense in $\cal R$.

Finally, the formula which follows directly from the
definition of universality shows us immediately that
the set of all universal matrices $\cal M$
is a $G_{\delta}$-set:

$${\cal M} =\cap_{k \in {\bf N}}\cap_{n \in {\bf N}}
\cap_{a \in A(r^n) \cap {\bf Q}^{n^2}} \cup_{m \in {\bf N}, m>n}\{r\in
{\cal R}:\max_{i=1,\dots n}
|r_{i,m}-a_i|< \frac {1}{k} \}.$$}
\end{proof}

Let us fix some infinite universal proper distance matrix $r$,
and provide the set of all natural numbers $\bf N$ with metric $r$.
Denote the completion of the space $({\bf N},r)$ with respect
to metric $r$ as a metric space $({\cal U}_r, \rho_r)$.
Evidently, it is a Polish space.

\begin{lemma}
{The distance matrix of any everywhere dense countable subset $\{u_i\}$
of the space ${\cal U}_r$  is a universal distance matrix.}
\end{lemma}

\begin{proof}
{Let us identify the set ${\bf N}$ with $\{x_i\} \subset {\cal U}_r$.
Then by definition $\rho(x_i,x_j)=r_{i,j}$
By definition the universality of $r$ means that for any $n$
the closure ($Cl$) in
${\bf R}^n$ of the set of vectors coincide with the set of the admissible
vectors of NW-corner of matrix $r$ of order $n$:
$Cl(\cup_{j>n}\{\{\rho(x_i,x_j)\}_{i=1}^n\}=A(p_n(r))$.
Because the set $\{u_i\}$ is also everywhere dense in $({\cal U}_r, \rho_r$
we can replace the previous set by the following:
$Cl(\cup_{j>n}\{\{\rho(x_i,u_j)\}_{i=1}^n\}=A(p_n(r))$. But because
$\{x_i\}$ is everywhere dense in $(U_r, \rho_r)$ we can change this on
$Cl(\cup_{j>n}\{\{\rho(u_i,u_j)\}_{i=1}^n\}=A(p_n(r'))$, where $r'$
is distance matrix for $\{u_i\}$.}
\end{proof}

We will see that $({\cal U}_r, \rho_r)$ is the so called Urysohn space
which is defined below, and the universality of matrix is necessary,
and sufficient condition to be a countable everywhere dense set
in Urysohn space.

\subsection{Urysohn universal space and universal matrices.}

Now we introduce the remarkable Urysohn space.
In his last papers \cite{U} Pavel Samuilovich Urysohn (1898-1924) gave
a concrete construction of the universal Polish space which is now called
 "Urysohn space". It was the answer on the question whihc was posed to him
by M.Freshet about universal Banach space. Later Banach and Mazur
have proved existence of the universal Banach spaces (f.e. $C([0,1]$),
but Urysohn's answer was more deep because his space is
homogeneuos in some sense. Actually Urysohn had proved several theorems
which we summarize as the following theorem:

\begin{theorem}(Urysohn \cite{U})

{A.There exists a Polish(=separable complete metric)
space $\cal U$ with the properties:

1)(Universality) For each Polish space $X$ there exists
the isometric embedding of $X$ to the space $\cal U$;

2)(Homogeneity) For each two isometric finite subsets $A=(a_1 \dots a_m)$
and $B=(b_1 \dots b_m)$ of $\cal U$ there exists isometry $J$
of the whole space $\cal U$ which brings $A$ to $B$: $JX=Y$ ;

B.(Uniqueness) Two Polish spaces with the previous properties 1) and
2) are isometric}.
\end{theorem}

The condition 2)  of the theorem could be strenghened: the finite
subsets possible to change to the compact subsets. So the group of isometry
of the space acts transitively on the isometric compacts.
But to enlarge the compact sets onto closed subsets is impossible.
Also we can in equivalent manner formulate this condition as
condition of the extension of the isometries from compacts to
the whole space, see below.

The main result of this section is the following
theorem, which includes the previous theorem.

\begin{theorem}
{1.The completion $({\cal U}_r, \rho_r)$ of the set of natural numbers
$({\bf N},r)$ with respect
to any universal proper distance matrix $r$ satisfies to the properties
1) and 2) of the above theorem and consequently is the Urysohn space.

2.(Uniqueness). For any two universal distance matrices $r$ and $r'$
the completions of the spaces $({\bf N},r)$ and  $({\bf N},r')$
are isometric.}
\end{theorem}

\begin{corollary}
 {The isometric type of the space $({\cal U}_r, \rho_r)$ does not depend
on the choice of universal matrix $r$.
The universality of the distance matrix is necessary
and sufficient condition to be the distance matrix
of a countable everywhere dense subset of the Urysohn space.}
\end{corollary}

The corollary follows from the Theorem 3 and Lemma 6.
The proof of the Theorem 3 which we give here repeats and simplifies
some arguments of Urysohn but he did not use infinite distance
matrices and the useful notion of the universal matrix.

\begin{proof}
{Suppose that matrix $r=\{r(i,j)\}_{i,j=1}^{\infty} \in {\cal R}$ 
be an arbitrary
universal proper distance matrix (it is convenient now
to write $r(i,j)$ instead of $r_{i,j}$)
and the space ${\cal U}_r$ is the completion of the countable
metric space $({\bf N},r)$, denote the corresponding metric on ${\cal U}_r$ as
$\rho_r$, but we will omit the index $r$.

1. First of all we will prove that the
metric space $(\cal U,\rho)$ is universal in the sense of property 1)
of the theorem 2, and then that it is  homogeneous in the sense of
property 2) of theorem 2.

Let $(Y,q)$ be an arbitrary Polish space. In order to prove
there is an isometric embedding of $(Y,q)$ into $(\cal U,\rho)$
it is enough to prove that there exists an isometric  embedding of
some countable everywhere dense set
$\{y_n\}_1^{\infty}$ of the space $(Y,q)$
to $({\cal U},\rho)$. This means that we must prove that for any infinite
proper distance
matrix $q=\{q(i,j)\} \in \cal R$ there exists some
countable set $\{u_i\} \subset \cal U$
with distance matrix equal to $q$. In it turn for this we need to construct a set
of the fundamental sequences in the space $({\bf N}, r)$, say,
$N_i=\{n^{(m)}_i\}_{m=i}^{\infty} \subset {\bf N}, i=1,2 \dots$
such that
$$\lim_{m \to {\infty}} r(n^{(m)}_i, n^{(m)}_j)=q(i,j), i,j=m,
\dots 
\mbox{and for all} \quad i,
\lim_{m,k \to \infty} r(n^{(m)}_i,n^{(k)}_i)=0.$$

 The convergence of each sequence $N_i=\{n^{(m)}_i\}_{m=i}^{\infty}$
in the space $({\cal U},\rho)$ when $m \to \infty$ to some point
$u_i\in {\cal U}, i=1,2 \dots $
follows from the fundamentality of this sequence e.g. from the
second equality above, and because of the first equality, the
distance matrix of the limit set $\{u_i\}$ coincided with the matrix $q$.
Now we we will construct the needed sequences
$\{N_i\}_{i=1}^{\infty} \subset \bf N$
by induction.
Choose arbitrarily a point $n^{(1)}_1 \in \bf N$,
and suppose that for given $m>1$
we already have defined the finite fragments
$L_k=\{n^{(k)}_i\}_{i=1}^k
\subset {\bf N}$
of the first $m$ sets
$\{N_i\}_{i=1}^m$,
for $k=1, 2 \dots m$ with property
$\max_{i,j=1,\dots k} |r(n^{(k)}_i, n^{(k)}_j)-q(i,j)|=\delta_k < 2^{-k},
k=1, \dots m$, and the sets $L_k$ mutually do not intersect.

Our construction of the set $L_{m+1}$ will depend only on the set $L_m$,
so we can for simplicity change the notations and renumber
$L_m$ as follow:
$n^{(m)}_i=i, i=1, \dots m$.

Now we will construct a new set $L_{m+1}=\{n^{(m+1)}_i\}_{i=1}^{m+1}
\subset \bf N$ with the needed properties in the following way. Consider
the finite metric space $(V,d)$ with $2m+1$ points
$y_1, \dots y_m; z_1, \dots z_m, z_{m+1}$
with the distances: $ d(y_i,y_j)= r_{i,j}, i,j=1 \dots m,
d(z_i,z_j)= q(i,j), i,j=1,\dots m+1; d(y_i,z_j)=q(i,j),
i=1 \dots m, j=1\dots m+1;
i \ne j,  \qquad
d(y_i, z_i)= \delta_m, i=1 \dots m.$ for some $\delta_m$.
It is easy to check that this is correct definition of the distances.
Denote the distance matrix of the space $(V,d)$
as $q_m$. Now apply corollary 2 ($\epsilon$-extension of isometry)
and enlarge the set $L_m=\{1, 2 \dots m\}$
with the new set $L_{m+1}$ with $m+1$ natural numbers
$\{n^{(m+1)}_i\}_{i=m+1}^{2m+1} \subset \bf N$
in such a way that the distance matrix of $L_{m+1}$ differ from
the NW-corner of the order $m+1$ of the matrix $q$ not more than
$\delta_m$
which is less than $2^{-(m+1)}$:
$\max_{i,j} |r(n^{(m+1)}_i, n^{(m+1)_j}) - q_m(i,j)|
=\delta_{m+1}<2^{-(m+1)}$
(remember that NW corners of order $m$ of matrices $q_m$ and
$r$ are coincided by construction).
We can see that for each $i$ the sequences
$\{n^{(m_i)}\}_{m=i}^{\infty}$
is fundamental and
$\lim_{m \to \infty} r(n^{(m)}_i,n^{(m)}_j)=q(i,j)$. Thus we have proved
that each Polish space can
be isometrically embedded into $({\cal U}_r,\rho_r)$.

We can essentially refine now the corollary 2 as follow.

\begin{corollary}(extension of isometry).
{The space $({\cal U}_r,\rho_r)$ has the following property: for any finite
set $A=\{a_i\}_{i=1}^n \in {\cal U}_r$ and distance matrix
$q$ of order $N, N>n$
with NW-corner of order $n$ which is equal to the matrix
$\{\rho(a_i,a_j)\}_{i,j=1}^n$ there exist points $a_{n+1} \dots a_N$
such that distance matrix of whole set $\{a_i\}_{i=1}^N$ is equal to 
the matrix $q$.}
\end{corollary}

The proof follows to the proof of Corollary 2 and uses the arguments
which we use above.

\begin{remark}
As we have mentioned there exist examples of universal but not homogeneous
Polish spaces (f.e. Banach space $C([0,1])$). The corollary above shows that
the main difference between such universal spaces and Urysohn space is the following:
we can isometrically embed in any universal space a given metric space;
but in the case of Urysohn space we can do more: the points of the  of embedded
metric space have given distances from the points of a fixed finite (or even compact)
set.
\end{remark}
 Let us continue the proof of the Theorem 3.

2. In order to prove homogeneity let us fix two finite $n$-point sets
$A=\{a_i\}_{i=1}^n$  and $B=\{b_i\}_{i=1}^n$ of  $({\cal U}_r,\rho_r)$
and construct two isometric ordered countable subsets  $C$ and $D$
each of which is everywhere dense in ${\cal U}$ and $C$ begins with
$A$ and $D$ begins with $B$.
The method of constrution is well-known and called "back and forth".
First of all we fix some countable everywhere dense subset $F$ in
$({\cal U}_r,\rho_r), F \cap A =F \cap B=\emptyset $, and represent
it as increasing union:
$F=\cup F_n$. Put $C_1=A \cup F_1$, and find a set  $D_1=B \cup F'_1$
so that the isometry of $A$ and $B$ extends to $F_1$ and $F'_1$.
Thus $D_1$ is isometric to $C_1$. This is possible to define because
of Corollary 4 (extension of isometry).
Then, choose $D_2=D_1 \cup F_2$ and $C_2=C_1 \cup F'_2$
and again extend the isometry  from the part on which it was defined before
to whole set. So we construct an isometry between $D_2$ and $C_2$ and so on.
The alternating process gives us two everywhere dense isometrical sets
$\cup C_i$ and $\cup D_i$ and the isometry between them extends isometry of $A$
and $B$.

3.Uniqueness. Let $r$ and $r'$ be two universal proper distance matrices and the spaces
$({\cal U}_r, \rho_r)$ and $({\cal U}_r', \rho_r')$  their completions.
We will construct repectively in the spaces two countable everywhere dense sets $F_1$
and $F_2$ so isometry between them will extend to the whole space.
Denote by $\{x_i\}$ and $\{u_i\}$ everywhere dense subsets of
$({\cal U}_r, \rho_r)$ and $({\cal U}_r', \rho_r')$ to which are generated
respectively by matrices $r$ and $r'$.
Now we repeat the same arguments as in the proof of the first
part of the theorem. We start with finite number of the points $\{x_i\}_{i=1}^{n_1}$
in $({\cal U}_r, \rho_r)$ and append to them the set of points
$\{u'_i\}_{i=1}^{m_1} \subset {\cal U}_r$ with the same distance matrix
as the distance matrix of the set of points $\{u_i\}_{i=1}^{m_1}$;
this is possible because of universality of the $({\cal U}_r, \rho_r)$ (property 1)
which had been already proved). Then append to the set $\{u_i\}_{i=1}^{m_2}, (m_2 > m_1)$
in the space ${\cal U}_r'$ the set of points $\{x'_i\}_{i=1}^{n_2}, (n_2 > n_1)$
in such a way that the distance matrix of the subset  $\{u_i\}_{i=1}^{m_1} \cup
\{x'_i\}_{i=1}^{n_1}$ of the set $\{u_i\}_{i=1}^{m_2} \cup
\{x'_i\}_{i=1}^{n_2}$  coincides with the distance matrix of the
set $\{u'_i\}_{i=1}^{m_1} \cup \{x_i\}_{i=1}^{n_1}$ etc.
continuing this process ad infinity as the result
of this construction we obtain two sets -
the first is $\{x_i\}_{i=1}^{\infty} \cup \{u'_i\}_{i=1}^{\infty} \subset {\cal U}_r$
and the second is
$\{u_i\}_{i=1}^{\infty} \cup \{x'_i\}_{i=1}^{\infty} \subset {\cal U}_r'$  -
which are everywhere dense in their spaces and are isometric. Thus
we have concluded the proof of the theorem.}
\end{proof}

Theorems 1 and 3 give us the following  the remarkable fact:
\begin{corollary}
{A generic ("typical") distance
matrix is a universal matrix, and consequently a generic
Polish space (in the sense our model of the cone $\cal R$)
is the Urysohn space $\cal U$.}
\end{corollary}

In his paper P.Urysohn gave an example of a countable space with rational
distance matrix (indeed that was universal incomplete metric space
over rationals $\bf Q$).
Our method of construction is more general: we construct the
universal matrix based on the geometry of the cone $\cal R$
and allows to give necessary and sufficient condition on
the distance matrix of any countable everywhere dense sets.
In the section  5 we apply it to the construction of Urysohn space
in probabilistic terms. We will give also the measure theoretic
versions of the universality of Urysohn space and prove some
facts about metric spaces with measure.

Urysohn also pointed out that there exist universal metric space of
the given diameter (say, equal to 1).
If we define in the same spirit the notion of universal matrix
with entries from interval $[0,1]$,
we obtain the universal metric space  of diameter 1
and the assertions of all theorems of this
paragraph take place for that space.

\section{Matrix distribution as complete invariant of the metric triple
and its characterization.}
\subsection{Matrix distribution and Uniqueness Theorem}

Now we begin to consider the metric spaces with measure and the
random metrics on the natural numbers.

Suppose $(X,\rho,\mu)$ is a Polish spaces with metric $\rho$ and with
borel probability measure $\mu$. We will call {\it metric triple}
(In \cite{G} the author used term ``mm-space'' another term is ``probability metric space'').
 Two triples  $(X_1,\rho_1, \mu_1)$ and
$(X_2, \rho_2, \mu_2)$ are {\it isomorphic} if there exists isometry
$V$ which preserve the measure:

$$\rho_2(Vx,Vy)=\rho_2(x,y),  V \mu_1=\mu_2.$$
As it was mentioned the classification of the Polish space (non compact)
is non smooth problem.
Surprisingly the classification of metric triple is "smooth"
and has a good answer which connects with the action of the group
$S_{\infty}$  and $S_{\infty}$-invariant measures on the cone $\cal R$.

For the metric triple $T=(X,\rho, \mu)$ define the infinite product
with the Bernoulli measure
$(X^{\bf N},\mu^{\bf N})$
and the map $F: X^{\bf N} \to {\cal R}$ as follows

$$ F_T(\{x_i \}_{i=1}^\infty) =\{\rho(x_i,x_j)\}_{i,j=1}^\infty \in {\cal R}$$
The $F_T$-image of the measure $\mu^{\bf N}$ which we denoted as $D_T$ will be called
{\it matrix distribution of the triple $T$}:
$F_T \mu^{\infty} \equiv D_T$.

The group $S_{\infty}$ of all finite permutations of the natural numbers
(infinite symmetric group $S_{\bf N}$) acts on the $\bf M_N(R)$ as well
as on the cone $\cal R$ of the distance
matrices as the group of simultaneous permutations of rows and
columns of matrix.
\begin{lemma}
{The measure $D_T$ is a Borel probability measure on  $\cal R$
which is invariant and ergodic with respect to the action of
infinite symmetric group, and invariant and ergodic
with respect to simultaneous shift
in vertical and horizontal directions (shortly NW-shift):
$(NW(r))_{i,j}=r_{i+1,j+1}; i,j=1,2 \dots$).}
\end{lemma}

\begin{proof}
{All facts follow from the same properties of the measure $\mu^{\infty}$,
which is invariant under the shift and permutations of the coordinates,
and because map $F_T$ commutes with action of the shifts and permutations.}
\end{proof}

Let us call a measure on the metric space non-degenerated
if there are no nonempty open sets of zero measure.
\begin{theorem}
{Two metric triples $T_1=(X_1,\rho_1, \mu_1)$ and $T_2=(X_2,\rho_2,\mu_2)$ with
non-degenerated measures  are equivalent iff its matrix distributions coincide:
$D_{T_1}=D_{T_2}$ as measures on the cone $\cal R$.}
\end{theorem}

\begin{proof}
{The necessity of the coincidence of the matrix distributions is evident:
if there exists an isometry $V:X_1 \to X_2$ between $T_1$ and $T_2$
which preserves measures then the infinite power $V^{\infty}$ preserves
the Bernoulli measures: $V^{\infty}(\mu_1^{\infty})=
\mu_2^{\infty}$ and because of equality
$F_{T_2}X_2^{\infty}=F_{T_2}(V^{\infty}X_1^{\infty}$),
the images of these measures are the same: $D_{T_2}=D_{T_1}$.
Suppose now that  $D_{T_2}=D_{T_1}=D$. Then $D$-almost all distance matrices
$r$ are the images under the maps $F_{T_1}$ and $F_{T_2}$, say
$r_{i,j}=\rho_1(x_i,x_j)=\rho_2(y_i,y_j)$ but this means that the identification
of $x_i \in X_1$ and  $y_i\in X_2$ for all $i$ is an isometry $V$
between these countable
sets. The crucial point of the arguments: by the ergodic (with respect to NW-shift) theorem
$\mu_1$-almost all sequences $\{x_i\}$ and $\mu_2$-almost all sequences
$\{y_i\}$ are uniformly distributed on $X_1$ and $X_2$ respectively.
This means that the $\mu_1$ measure of each ball
$B^l(x_i) \equiv \{z\in X_1:\rho_1(x_i,z)<l\}$
is equal to $$\mu_1(B^l(x_i))
=\lim_{n \to \infty}\frac{1}{n} \sum_{k=1}^n 1_{[0,l]}(\rho_1(x_k,x_i)). $$
But because of the isometry $V$ ($r_{i,j}=\rho_1(x_i,x_j)=\rho_2(y_i,y_j)$
- see above) the same quantity is a $\mu_2$-measure of the ball:
$B^l(y_i) \equiv \{u\in X_2:\rho_2(y_i,u)<l\}$
 $$\mu_2(B^l(y_i))
=\lim_{n \to \infty}\frac{1}{n} \sum_{k=1}^n 1_{[0,l]}(\rho_2(y_k,y_i))= \mu_1(B^l(x_i)). $$
Finally, both measures are non-degenerated, consequently each of the sequences $\{x_i\}$
and $\{y_i\}$ is everywhere dense in its own space. Because
both measures are Borel it is enough to conclude their coincidence
if we establish that the measures of the all such balls are the same.}
\end{proof}

\begin{corollary}
{Matrix distribution is complete invariant of the equivalence classes
(up-to isometries which preserve the measure) of the of metric triples
with non-degenerated measures.}
\end{corollary}

  We can call this theorem the ``Uniqueness Theorem'' because it
asserts the uniqueness up-to isomorphism of the metric triple with
the given matrix distribution.  Firstly this theorem
as the ``Reconstruction Theorem'' in another formulation
has been proved in the book \cite{G} pp.117-123 by Gromov.
He formulated it in the terms of finite dimensional distributions
of what we called matrix distribution and proved it using
analytical method. He asked me in 1997 about this theorem
and I suggested  the proof which is written here (see also in \cite{V1})
and which he had quoted (pp.122-123) in the book.
Gromov had invited the readers to compare two proofs, one of which is
rather analytical (Weierstrass approximations) and another (above)
in fact uses only the reference to the
ergodic theorem. The explanations of this
difference is the same as in all applications of the ergodic theorem -
it helps us to replace the methods of space approximation by operations
with infinite (limit) orbits. In our case the consideration of infinite
matrices and cone $\cal R$ with invariant measures gives a possibility
to reduce the problem to the investigation of ergodic action
of infinite groups. For example the uniformicity of the distribution
of the sequence has no meaning for finite but very useful for infinite
sequences.
In \cite {V2,V3} we use a more general technique which is also based on the ergodic
methods in order to prove the analog of uniqueness theorem for the classification
of arbitrary measurable functions of two variables
(in the case above this was a metric as a function of two arguments).

\subsection{Properties of matrix distributions and Existence Theorem}

The matrix distribution of a nondegenerated metric triple
$T=(X,\rho, \mu)$
is by definition the measure  $D_T$ on the cone $\cal R$.
Clearly it can be cosidered as a random (semi)metric on 
the set of natural numbers. In this section we will characterize
those random metrics (or those measures on $\cal R$)
which could be a matrix distribution, in other words those 
distributions on the cone $\cal R$ which can appear
as a random distance matrices for independent sequences of points
$\{x_i\}$ of the metric space $(X,\rho)$ which are distributed with
measure $\mu$ on that space. To characterize this set is necenssary in
orfer to cliam that the classification problem is indeed smooth and
the set of the invariants has an explicit description.
We will show that the set of matrix distributions is the borel set 
in the space of probability measures on the cone $\cal R$.
 
As we mentioned (Lemma 7) any measure $D_T$ must be invariant and
ergodic with respect to action of infinite symmetric group and
to NW-shift. But this is not sufficient and below one can find 
the necessary and sufficient conditions for that (see also \cite{V3}).
But will start from the counterexamples.

{\bf Examples.}

1.A trivial example of an invariant ergodic measure which is not
matrix distribution
is the following. Denote $r^0$ a distance matrix:
$r^0_{i,j}=\delta (i-j)$ (where $\delta (n)=1$ if $n=0$ and $=0$
otherwise);  this
is nothing than distance matrix of the countable set such that
the distance between two different
points is equal to $1$. Let a measure $\mu^0$ be a delta-measure
at the matrix $r^0$.
Clearly $\mu^0$ is invariant, ergodic and does not correspond
to any metric triple.

2. An example of the general type is the following. Firstly note
that each symmetric matrix with zeros on the principal
diagonal and with entries $r_{i,j}$ from the interval $[1/2,1]$
when $i \ne j$ is a proper distance matrix; indeed in this case
the triangle inequality is valid for each three numbers.
For each probability measure $m$ with the support
on  $[1/2,1]$ which is not just a single atomic measure consider
the product measure $m^{\infty}$ with the factor $m$ on
${\bf M}^s_{\bf N}(\bf R)$ (it means that all entries upper diagonal
are independent and identically distributed). Consequently this measure
is concentrated on the cone $\cal R$. This measure is evidently
invariant under permutations and under NW-shift as well as
is ergodic measure with respect to those transformation.
In the same time this continuous measure is not a matrix distribution
for any metric triple. Indeed, with
$m^{\infty}$-probability equal to one all distance matrices
define the discrete topology on the set of natural numbers $\bf N$
because of the absence of nontrivial fundamental sequences in  $\bf N$,
and consequently the completion of $\bf N$ is $\bf N$, thus
a matrix distribution cannot be a continuous measure, but our
product measure $m^{\infty}$ is continious one.

The explanation of those effects will be clear from the proof
of the next theorem which gives one of the  characterizations
o matrix distribution.

\begin{theorem}(Existence of metric triple with given matrix distribution)

{Let $D$ be a probability measure on the cone $\cal R$,which is invariant
and ergodic with respect to action of infinite symmetric group
(=group of all finite permutations of the naturals).

1)The following condition is necessary and sufficient for $D$ to be a matrix distribution
for some metric triple $T=(X,\rho, \mu)$ or $D=D_T$:

for each $\epsilon>0$ there exists integer $N$ such that
\begin{equation}
D\{r=\{r_{i,j}\} \in {\cal R}: \lim_{n \to \infty} \frac{|\{j:1\leq\ j \leq n,
\min_{1 \leq i \leq N}
 r_{i,j} <\epsilon \}|}{n}> 1- \epsilon \}>1-\epsilon.
\end{equation}

2)The following stronger condition is necessary and sufficient for $D$
to be a matrix distribution for some metric triple
$T=(X, \rho, \mu)$ with compact metric space $(X, \rho)$:

for each $\epsilon >0$ there exists integer $N$ such that

\begin{equation}
D\{r=\{r_{i,j}\} \in {\cal R}: \mbox{for all} \, j>N,
\min_{1\leq i \leq N} r_{i,j} <\epsilon\}>1-\epsilon,
\end{equation}
}
\end{theorem}

\begin{proof}
{A.Necessity. In the case of compact space the necessity is evident:
the condition (5) expresses the fact that sufficiently long
sequences of independent (with respect to the nondegenerated $\mu$)
points being uniformly distributed with respect to $\mu$
contain an $\epsilon$-net of the space.
The necessity of conditions (4) follows automatically from
well-known property of the borel probability measures
on the complete separable metric space:
namely a set of full measure is sigma-compact
(so called ``regularity of the measure''),
consequently for each $\epsilon$ there exists
a compact of measure $>1-\epsilon$. Indeed, because of countably additivity
of our measure  for any $\epsilon >0$ the exists finite number of the points
such that the measure of the union of $\epsilon$-balls with the centers at those
points is greater than $1-\epsilon$ and using a ergodic theorem
we can assert that the condition in the brackets in (4) valids
for matrix distance from the set of measure more than  $1-\epsilon$.

B.Sufficiency Suppose now that we have a invariant and ergodic measure $D$ on $\cal R$ with
condition (4). The plan of the proof is the following: we express
all the properties of the measure $D$ in terms of "typical" distance
matrix  $r$ and then we will construct
a metric space with measure (metric triple) using only one "typical"
distance matrix $r$.
 Invariance of $D$ under the group $S_\infty$ (simultaneous permutations
of the rows and columns) leads to the invariance of the restrictions
of the measure $D$ on the submatrix $\{r_{i,j}: i=1,2 \dots n, j=1,2\dots \}$
with respect to the shift $j \to j+1$ for any $n$.
Using ergodic theorem for this shift (which is not ergodic!)
we can find the set $F \subset \cal R$  of full $D$-measure
of such distance matrices $r=\{r_{i,j}\}$ for which the following limits exist for
any natural numbers $k$ and positive real numbers $\{h_i\} \, i=1,2 \dots k$:
\begin{equation}
   \lim_{n \to \infty}\frac{1}{n}\sum_{j=1}^n \prod_{i=1}^k 1_{[0,h_i]}(r_{i,j})
\equiv \mu^{h_1, \dots h_k}
\end{equation}

  Now let us use the invariance of the measure $D$ under the action
of symmetric group $S_\infty$.
By the ergodic theorem (more exactly by the martingale theorem) for the
action of $S_\infty$
as locally finite group we can assert that for almost
all $r$ and fixed Borel set $B \in {\cal R}_n$  the following limits exist
\begin{equation}
\lim_{N\to \infty}1/N!\sum_{g \in S_N}1_B (g(r^{(n)})) \equiv \Lambda^{(n)}_r(B)
\end{equation}
where $ g(r^{(n)})=\{r_{g(i),g(j}\}_{i,j=1}^n$, $g$ is a permutation i.e. element
of $S_N$, which permutes the first $N$ naturals numbers,  $1_B$ is a characteristic
function of the Borel set $B \subset{\cal R}_n $; the measure
 $\Lambda^{(n)_r}(.)$ in ${\cal R}_n$ is called  the empirical distribution
of matrix $r \in{\cal R}_n$.
These empirical distributions as a family of measures on the cones ${\cal R}_n$
are concordant with respect to the projections $p_n$ (see section 2) and consequently define
an $S_\infty$-invariant measure on $\cal R$. Our assumption about the matrix $r$
is that this measure coincides with the initial measure $D$; it is possible to assume
this because of ergodicity of action of $S_\infty$. If we choose a countable basis
of the Borel sets $\{B_i^n\}_{i=1}^{\infty}$ in ${\cal R}_n; n=1,2 \dots $
then for $D$-almost all $r$ the existence is valid for all $B_i^n, i,n=1,2 \dots $.

Finally let us restate the condition (4) in terms of the distance matrix $r$.
It follows from (4) that the for $D$-almost all $r$ the following is true:
for each $k$ there exist integer $N$ such that

\begin{equation}
 \lim_{n \to \infty} \frac{|\{j:1\leq\ j \leq n,
\min_{1 \leq i \leq N}
 r_{i,j} < k^{-1} \}|}{n}> 1- k^{-1}.
\end{equation}

Let us fix  one such distance matrix $r=\{r_{i,j}\}$ which satisfies to the
conditions (6-8) and consider it as a metric on the set of natural numbers.
Denote by  $X_r$ of the completion of the metric space $({\bf N},r)$,
denote the metric in this completion by $\rho_r\equiv \rho$,
and the natural numbers as a dense countable
set in this completion by $X_r$ by $x_1, x_2, \dots $.
Denote by $B^h(x)$  the ball of the
radius $h$ with the center at the point $x$ in the space $X_r$ and let
$\cal A$ be the {\it algebra of subsets} of $X_r$ generated
by all the balls with the center at the  points $x_i, i=1,2 \dots $
and arbitrary radius.
Let by definition the measure $\mu_r$ of the finite
intersections of the balls be as follows:

\begin{equation}
\mu_r(\cap_{i=1}^k B^{h_i}(x_i))=  \mu^{h_1, \dots h_k}
\end{equation}

It is easy to check that this equality correctly defines nonnegative
finitely additive normalized measure $\mu_r$ on the {\it algebra} $\cal A$
of the sets generated by the mentioned balls,  but in general it is NOT sigma-additive
and consequently can not be extended to sigma-algebra of all Borel set in $X_r$
as true probability measure. This is just the case in our
counterexamples above: we have had a countable space
and the definition above gave a measure which takes value zero
on each finite set but equal to 1 on the whole space.
\footnote{In a sense we are in the situation of the
classical Kolmogorov's theorem about extension of the
measures and its generalizations: the measure is defined
on the algebra of the cylindric sets and after the test
on countable additivity  we can extend a measure on sigma-algebras.
This is possible for each measures in the linear space $R^{\infty}$
(Kolmogoroff's theorem), but not possible in general in other spaces.
In our case the measures are defined on the algebra generated with balls
and condition (4) guarantees the countable additivity; in our
cases also for some spaces countable additivity takes place
automatically.}

Now we will use the condition (4) in the form (8) for $r$.
Choose  $\epsilon > 0$,
áondition (8) allows to find for each $k$ a finite
union of the balls, say, $C$ in $X_r$ of measure more than $1-\epsilon$.
Normalize our measure on  $C$ to $1$, denote it as $\bar \mu_r$.
Using induction on $k$ we can construct a set of balls with radius
$2^{-k}$ such that the intersection of union of $C_k$ and $C$ has
the $\bar \mu_ r$-measure more than $1-2^{-k}\epsilon$; $k=1,2 \dots$.
This means that the intersection of all these
sets $ C \cap (\cap_k C_k)$ has $\bar \mu_r$-measure more
than $1-2\epsilon$ and is a totally bounded set
(i.e.has an $\epsilon$-net for all
$\epsilon$); because of completeness of $X_r$ this intersection is a compact.
But any  finite additive measure which is defined on an algebra of the sets dense in the
sigma algebra of the Borel sets in the compact is countable additive. So we have found
a compact $C$ in $X_r$ whose $\mu$-measure is not
less than $1-3\epsilon$. Because
$\epsilon$ is arbitrary we have constructed atrue probability
measure $\mu$ in $X_r$  with sigma-compact support. If we use instead of
conditions (4) and (8) the condition (5) and its individualization for $r$
we obtain along the same construction a compact of full measure in $X_r$.

We have constructed a metric triple  $T_r=(X_r,\rho_r,\mu_r)$
where the measure $\mu$ is probability measure on the
Polish space $(X_r,\rho_r)$ with full support
and with distinguished dense countable subset $\{x_i\}$ which is
{\it uniformly distributed (with respect to measure} $\mu_r$), and
also satisfes  the condition (7). The final part of the proof consists
in the verification of the fact that matrix distribution $D_{T_r}$  of the metric
triple $T_r=(X_r, \rho_r, \mu_r)$ and initial measure $D$
are equal as the measures on the cone $\cal R$.
We formulate this as a Lemma which is useful in more general situations.
This completes the proof of the theorem.}
\end{proof}

\begin{lemma}
{Suppose $r \in \cal R$ is a matrix for which all the limits (6)
exist and equation (7),(8) is also valid.
Construct the metric triple
$T_r=(X_r, \rho_r, \mu_r)$ using
the equations (6) and (8).
Then the matrix distribution $D_{T_r}$ of this triple
is equal to the $S_{\infty}$-invariant measure which is generated
from the  matrix $r$ by formula (7).}
\end{lemma}

\begin{proof}
{For the proof we must check
the coincideness of the finite dimensional distributions of both measures.
Let us illustrate this for the case of the distribution of the element
$r_{1,2}, (n=2)$;  for general $n$ the verification is similar.
$$
\begin{array}{l}
\int_{X_r}\int_{X_r} {\bf 1}_B(\rho_r(u,z))d\mu_r(u)d\mu_r(z)=
\lim_{n \to \infty}n^{-1} \sum_{i=1}^n \lim_{m \to \infty}m^{-1}\sum_{j=1}^{m}
 {\bf 1}_B(r_{i,j})=\\
\lim_{n \to \infty}n^{-2}\sum_{i,j=1}^n {\bf 1}_B(r_{i,j})=
\lim_{N \to \infty}({N!})^{-1}\sum_{g \in S_N}{\bf 1}_B(r_{g(1),g(2)})=
\Lambda^{(2)}_r(B).
\end{array}.
$$
 Here $B \subset {\bf R}_+$; the last equality follows from (7);
the equalities above used the
uniformity of the distribution of the sequence $\{x_i\}$
in the space $(X_r, \mu_r)$.
This concludes the proof of the theorem, because by the condition (7)
the $S_{\infty}$-invariant measure on $\cal R$ which is generated by
matrix $r$ is just the measure $D$.}
\end{proof}

\begin{remark}
{1.The structure of the conditions on the measure in the thoerem shows that
the set of matric distributions is indeed a borel set in the space of all borel
probability measures on the cone $\cal R$.

2.The condition (4) could be replaced by another condition from the paper \cite{V3}
(simplicity of $S_{\infty}$-invariant measure). That condition guarantees
that measure $D$ appeared from {\it some} measurable function of two variables
as matrix distribution which is sufficient for our goals.}
\end{remark}

\subsection{The space of measure-theoretical metric triples.}

We can extend the notion of the space of metric spaces (see section 1)
and introduce a similar space for the metric triples. Instead of the ordinary
point of view where one considers the set of all Borel measures on the given topological
space, we in opposite, consider the set of all {\it measurable (semi)metrics}
on a fixed Lebesgue space with continuous measure. (see \cite{V4}, par.6).

  Let $(X,\mu)$ be a Lebesgue space with measure $\mu$ finite or
sigma-finite (say, interval [0,1] with Lebesgue measure or natural numbers
with  the uniform mreasure), and $S_{\mu}(X)$- the space of all classes
$mod 0$ of measurable functions;  define ${\cal R}^c \subset S_{\mu}(X)$
as a  cone of measurable metrics e.g. the cone of the classes $mod 0$
of symmetric measurable functions $\rho :(X\times X, \mu \times \mu)) \to {\bf R}_+$
with the triangle inequality:

$$\rho(x,y)+\rho(y,z) \geq \rho(x,z) \mbox {for} \qquad  (\mu \times \mu \times \mu)-
\mbox {almost all} \qquad  (x,y,z)\in (X \times X \times X).$$
It is natural to assume that $\mu \times \mu\{(x,y): \rho(x,y)=0\}=0$.

Remark that $\rho$ is not individual function but the class of $mod 0$
equivalent functions, so it is not
evident a priori that such $\rho$ defines the structure
of metric space on $X$ in the usual sense.

If measure  $\mu)$ is discrete one then the cone ${\cal R}^c$ is
the cone of ordinary (semi)metric on the finite or countable set
e.g. they coinside with  $\cal R$ or  ${\cal R}_n$ (see section 1),
Thus the cone ${\cal R}^c$ is a continuous generalization
of the cone $\cal R$ where instead of the set of natural numbers $\bf N$ with counting measure
we consider the space $(X,\mu)$ with continuous measure.

Suppose now that measure $\mu)$ is finite and continuous and
$\rho \in {\cal R}^c$ is a {\it pure} function
(see \cite{V3})
\footnote{A measurable function $f(x,y)$ of
two variables on the unit square with Lebesgue measure
calls pure if for almost all pairs $(x_1,x_2)$ (corresp. $(y_1,y_2)$)
the functions of one variable $f(x_1,.)$ and $f(x_2,.)$ (corresp. $f(.,y_1)$ and $f(.,y_2)$)
are not coniside everywhere. Evidently, a proper metric
on the measure space is pure function.},
and the measure $ D_{\rho}$ on the space $\bf M_{\infty}(\bf R)$
is a matrix distribution of the measurable function $\rho$.
(see definition in the previous section).
Using the ergodic theorem we can prove that
$D_{\rho}(r \in {\cal R}:r_{i,k}+r_{j,k}\geq r_{i,k})=1$  for each $i,j,k \in \bf N$
and consequently $D_{\rho}({\cal R})=1$. From this  using characterization of
matrix distributions from \cite{V4} we  conclude that the following
assertion is true:

\begin{lemma}
{The measure $ D_{\rho}$ concentrates on the cone $\cal R$ (e.g.$D_{\rho}({\cal R})=1$)
and is an ergodic $S_{\infty}$-invariant measure.
Consequently, each pure function $\rho \in {\cal R}^c$ on $(X \times X, \mu \times \mu)$
defines a true metrics $mod 0$ on the space $(X,\mu)$.}
\end{lemma}
\begin{corollary}
{The class of measurable (semi)metrics on the Lebesgue space with
continuous measures coinsides with the class of the metric triples
with finite continuous measures.}
\end{corollary}
This corollary shows that the language of the matrix distributions
which are concentrated on the cone  $\cal R$, is an invariant
manner to study (semi)metric triples. It seems that sometime
it is convinient to fix the measure space and to vary in measurable manner
the metrics instead of consideration of the fixed metric spaces with
various measures - the generality of the objects is the same. We have used
this way in \cite{V4}).

\newpage

\section{General classification of the measures on the cone
of distance matrices, examples.}
\subsection{Definitions}
Let us consider arbitrary measures on the cone $\cal R$,
or - arbitrary random metrics on the naturals and choose some notations.
Remark that the cone $\cal R$  with weak topology is metrizable
separable  space (e.g. is the Polish space if we fix a metric
which is compatible with weak topology).

{\bf Notation.}
Denote  by $\cal V$ the set of all probability Borel measures \footnote{We use later the
term ``measure'' in the meaning ``Borel probability measure''} on the cone $\cal R$
and endow it with weak topology,- this is the topology of inverse
limit of the sets of probability measures on the finite dimensional cones ${\cal R}_n$
with its usual weak topology. The convergence in this topology is
convergence on the cylindric sets with open bases.
All classes of measures which we define below are
the subsets of $\cal V$ with induced topology.
Remark that the set of non-degenerated
(=positive on the nonvoid open sets) measures
is of course everywhere dense $G_\delta$ set in $\cal V$.

 Let $\cal D$  be  the subset in $\cal V$ of the matrix distributions;
as we proved (Theorem 4,'orollary 7) this set is in the bijective
correspondense with the set of all classes of isomorphic all metric triples.
The constructive description of $\cal D$ follows from the
existence theorem (section 4).

The subset $\cal P$ of  $\cal V$ is the set of measures which
are concentrated on the set of universal distance matrices:
$\nu \in \cal P$ iff $\nu({\cal M})=1$ - see section 3.

 The subset $\cal W$ of  $\cal V$ is the set of measures which
are concentrated on the set of weakly universal distance matrices
Both sets are convex (not closed) subspaces of the simplex of
all measures on the cone $\cal R$.
It is possible to give direct characterization of those measures
which analogous to the criteria of universality from Statement 1.

\begin{statement}
{Let us define for each natural $n$ a partition of the cone $\cal R$
on the subsets ${\cal R}^n(q), q \in {\cal R}_n$ (see denotation
before Statement 1 in the section 3). Measure $\mu$ belongs to
the cone $\cal P$, (e.g. measure is concentrated on the set of
universal matrices iff for each $n$ and for almost all elements
of those partitions (or for almost all finite distance matrix $q$)
the support of {\it conditional} measure is whole set ${\cal R}^n(q)$
(conditional mesusre are not degenerated).
Belonging of the measure $\mu$ to the set $\cal W \subset V$
is equivalent to the following: support of measure $\mu$ is
whole cone $\cal R$.}
\end{statement}

The set ${\cal Q} \subset \cal D$ consists of the measures
which corresponds to the metric triples $T=({\cal U}, \rho, \mu)$
in which $\cal U$ is Urysohn space.

Finally denote by $\cal H$ the set of measures $\mu$ in  $\cal V$
which have the following property: $\mu$-almost all distance
matrices generate the isometric metric spaces. A measure
$\mu \in \cal H$ generates a {\it random everywhere dense
sequence of points} on the given space.
From this point of view the elemetns of $\cal D$ induced
{\it a random everywhere dense subset of the special type}, namely
infinite independent sampling of the points of the given metric
triple; and the elements of $\cal Q$ induced a random independent
sampling of the points in the Urysohn space with nondegenerated measure.

We have embeddings:

$$  {\cal V} \supset {\cal H} \supset {\cal D} \supset {\cal Q} \subset
{\cal P} \subset {\cal H} \subset {\cal V};
\qquad {\cal Q}= {\cal P}\cap{\cal D} $$

The Theorem 4 shows that each measure from the set $\cal D$
defines a class of isomorphic metric triples and in particular
the set $\cal Q$ is the set of classes isometric nondegenerated
measures on the Urysohn space  (e.g. the orbit of the
group all isometries of the Urysohn space on the set of
nondegenrated measures on that spaces.

\begin{theorem}
{1.The subset ${\cal P} \subset \cal V$
is an everywhere dense $G_\delta$ subset in $\cal V$.
This means that for generic measures $\nu$ on $\cal R$
$\nu$-almost all distance matrices
are universal, and consequently
                                                                                                                               tance matrix  $r$ defines
a metric on the naturals such that completion of
naturals with respect to this
metric is the Urysohn space.

2.The subset ${\cal Q} \subset \cal D$
is an everywhere dense $G_\delta$ in $\cal D$.
This means that the generic metric triple $T=(X, \rho, \mu)$
has Urysohn space as the space $(X,\rho)$.}
\end{theorem}

\begin{proof}
{The first claim follows from Theorem 1 which states in particular
that the set of universal distance matrices is a $G_\delta$
in $\cal R$, and from a general fact we can deduce that the set of all measures
on separable metrizable space such that some fixed everywhere dense $G_\delta$
(in our case - $\cal M$)
has measure 1, is in its turn itself an everywhere dense $G_\delta$
in the space of all measures in the weak topology.
The second claim follows from the fact that the intersection
of $G_\delta$-set with any subspace in a Polish space is $G_\delta$ in
induced topology.}
\end{proof}

\subsection{Examples of the measures which are concentrated
on the universal matrices.}

Now we can give a probabilistic (markov) construction of the measures
on the cone $\cal R$ and in particular to represent the examples
of the measures from the set $\cal P \subset \cal H \subset \cal V$.
This gives a new proof of the existence
of the Urysohn space. In fact we use the arguments from the section
3 but in the probabilistic interpretation.  Also the method gives
a concrete illustration how to construct a random metric space.

Let $\gamma$ be an arbitrary continuous measure on half-line ${\bf R}^1_+$
with full support - for example Gaussian measure on the half-line.
We will define inductively the measure $\nu$ on the cone of distance matrices
$\cal R$ by construction of its finite dimensional projection on
the cones ${\cal R}_n^0$ or in other words - joint distributions
of the elemetns of random distance matrices.

The distribution of the element $r_{1,2}$ of the random matrix
is distribution $\gamma$. So we have defined the measure on ${\cal R}_2$, denote it as $\nu_2$.
  Suppose we have already defined the joint distribution of the
entries $\{r_{i,j}\}_{i,j=1}^n$, which means that we have defined
a measure $\nu_n$ on ${\cal R}_n^0$.
 By Lemma 4 the  cone  ${\cal R}_{n+1} $ is a fibration over
cone ${\cal R}_n$ with the fibers $A(r)$ over matrix $r\in {\cal R}_n$
We use only the structure of this fibration:
projection ${\cal R}_n\stackrel{p_{n+1}}{\longleftarrow}{\cal R}_{n+1}$ -
in order to define a measure on ${\cal R}_{n+1}$ with given projection.

So we need to define a conditional measure on $A(r)$ for all
$r \in \cal R$ which are measurably depend on $r$.
From probabilistic point of view this means
that we want to define transition probabilities from given
a distance matrix $r$ of order $n$ to the distance matrix $r^a$
(see section 2) of order $n+1$.

Let us recall the geometrical structure of the set of admissible
vectors $A(r)$. It is Minkowski sum: $$A(r)= M_r +\Delta_n,$$
(see 2.2) or as projection
of the direct product  $\pi:M_r \times \Delta_n \to  M_r +\Delta_n=A(r)$.
Consider product measure on
$M_r \times \Delta_n :\gamma_r=m_r \times \gamma$
where $m_r$ is for example normalized Lebesgue measure on the compact
polytope $M_r$ or another measure with full support on $M_r$,
with the conditions which we formulate below.

Let $\pi \gamma_r$ be its projection on $A(r)$.
 We define the {\it conditional measure} on $A(r)$ as $\pi \gamma_r$.
So, we have
$$ \mbox{Prob}(r^{da}|r)=\pi(m_r \times \gamma)(da). $$
The conditions on the measures $m_r$ are the following:
at each step of the construction
for each $N$ and $n>N$ the projection of the measure
$m_r, r \in {\cal R}_n$ to the set of admissible vectors $A(p_N(r))$
are uniformly positive on the open sets; this means that for any
open set $B \subset A(p_N(r))$ there exist
$\epsilon >0$ such that for any $n>N$ the value of projection of the measure
$m_r, r \in {\cal R}_n$ on the set $B$ more than $\epsilon$.
Thus we define a measure $L_n$ on ${\cal R}_{n+1}$.
By construction all these measures are concordant
and define the a measure $L$ on $\cal R$.
Denote this measure as
$L=L(\gamma, \{m_r:r \in {\cal R}_n, n=1,2 \dots \})$.

 A more intuitive and combinatorial variant of this description
is the following:
to the given $n$-point metric space we randomly add a $n+1$-th  point
choosing the vector of the distances between the new and the
previous points (admissable vector), with  the natural probability
which is positive on all open sets of admissible vectors.

\begin{theorem}
{The constructed measure $L$ belongs to the cone $\cal P$
which means that $L$ is concentrated on the set of universal matrices.
and therefore the completion of $({\bf N},r)$ is Urysohn space $L$-almost sure.}
\end{theorem}

\begin{proof}
{It is sufficient to check the conditions of the Statement 2.
By the theorem about convergence of the martingales the conditional
measure on almost all elements of the partitions -  ${\cal R}^n(q)$ -
is the limit of the conditional measures on the elements of the
partitions ${\cal R}^n(q) \cap {\cal R}_N$ when $N \to \infty$.
Thus the claim that almost all conditional measures are
nondegenerated is the consequence of the condition about
uniform positivity of the probability on the set of admissable vectors.}
\end{proof}

We can say that {\it a random countable metric space is an
everywhere dense subset of the Urysohn space}
or equivalently completion of the random countable metric space
with probability one is the the Urysohn space.
Here ``random'' means randomness respectively to that
natural procedure which was defined above and which is in a
sense very close to independence and has a very wide variations
which allow to define the measures on $\cal R$.

A much more complicated problem is to construct a measure on $\cal R$
from the set $\cal Q $ i.e. which is a matrix distribution for some
measure on $\cal U$. The properties of the measures on the Urysohn space
are very intriguing. But becuase we have no useful model for thisspace,
it is natural to use indirect way for the definition and studying
of such measures: to define matrix distribution as a measure on the
cone  $\cal R$ which belongs to the set ¢ $\cal Q$ which defined
an isometric class of the measures on $\cal U$. In its turn for this
we can take any measure from the set $\cal P$ and then to construct
the $S_{\infty}$-hull of its -$S_{\infty}$-invariant ergodic measure
The simplification is that we can omit the condition (4) from the
theorem which guarantees the fact that the measure  on $\cal R$
is matrix distribution:

\begin{statement}
{Each $S_{\infty}$-invariant ergodic measure from the set $\cal P$
is matrix distribution (belongs to $\cal Q$). This means that
$S_{\infty}$-hull of measures which we had constructed
above defines the isometry class of the measures on Urysohn space.}
\end{statement}
The proof is based on the criteria of simplicity of the measures
on the set of infinite matrices from \cite{V3}, and we will
discuss it elsewhere.

The probabilistic analysis on the distance matrices is useful
for integration over set of metric spaces in the spirit
of statistical physics. The measures which we had considered here
are interesting from the point of view of the modern theory
of random matrices. It is natural to study the spectra of the
random distance (symmetric) matrices. The simpliestcase is to
calculation of the joint distributions of distances between
independent random points of the homogeneous manifolds spaces
(spheres, for example), this is interesting and complicate
new problem.

Returning back to the Theorem 7, I must recall a very interesting
analogy with the old and simple theorem by  Erd{\"o}s-R\'{e}nyi
\cite{ER} about the random graphs. It asserts that with probability one
a random graph is the universal graph, see \cite{Ra,Ca}.
This is the simplest case of the the theory which we developed here
because each infinite graph defines the distance onthe set of vertices
and in our case ditance takes only two nonzero values - $1$ and $2$.
The random graph in the sense of \cite{ER} defines the measure on
the distance matrices such that all the entries of matrix are indepedent
and have unifirm distribution on $\{1,2\}$ (of more general).
All the matrices belong to the cone $\cal R$. In the same time almost
all matrices are universal in our sense (see paragraph 3)
if we consider only two values $\{1,2\}$ of distances
(instead of values from ${\bf R}_+$).

\bigskip

\hbox{A.Vershik.}
\hbox{St.Petersburg Math.Inst. of Russian Acad Sci.}
\hbox{vershik@pdmi.ras.ru}

\end{document}